\documentclass[11pt]{amsart}
\usepackage{palatino}
\usepackage[all]{xy,xypic}

\usepackage{amsfonts,amsmath,oldgerm,amssymb,amscd}

\newcommand{\ra}{\rightarrow}		
\newcommand{\lra}{\longrightarrow}
\newcommand{\by}[1]{\stackrel{#1}{\ra}}

\newcommand{\remove}[1]{}

\newcommand{\surj}{\ra\!\!\!\ra}	

\newcommand{\ol}{\overline}		
\newcommand{\wt}{\widetilde}
\newcommand{\iso}{\by \sim}

\newtheorem{theorem}{Theorem}[section]
\newtheorem{proposition}[theorem]{Proposition}
\newtheorem{lemma}[theorem]{Lemma}
\newtheorem{definition}[theorem]{Definition}
\newtheorem{corollary}[theorem]{Corollary}

\newcommand{\ga}{\alpha}	\newcommand{\gb}{\beta}
		\newcommand{\gd}{\delta}
	
	\newcommand{\gl}{\lambda}
	
		\newcommand{\gs}{\sigma}

\newcommand{\ul}{\underline}

\newcommand{\BC}{\mbox{$\mathbb C$}}

	\newcommand{\BP}{\mbox{$\mathbb P$}}
\newcommand{\bq}{\mbox{$\mathbb Q$}}	\newcommand{\br}{\mbox{$\mathbb R$}}

	\newcommand{\bz}{\mbox{$\mathbb Z$}}

\newcommand{\hh}{\text{ht}}

\newcommand{\dd}{\text{dim}}

\newcommand{\sur}{\twoheadrightarrow}
\newcommand{\bp}{\begin{proposition}}
\newcommand{\ep}{\end{proposition}}
\newcommand{\bl}{\begin{lemma}}
\newcommand{\el}{\end{lemma}}
\newcommand{\bt}{\begin{theorem}}
\newcommand{\et}{\end{theorem}}
\newcommand{\bc}{\begin{corollary}}
\newcommand{\ec}{\end{corollary}}
\newcommand{\bd}{\begin{definition}}
\newcommand{\ed}{\end{definition}}
\newtheorem{conjecture}[theorem]{Conjecture}

\def\rmk{\refstepcounter{theorem}\paragraph{{\bf Remark} \thetheorem}}
\def\proof{\paragraph{Proof}}

\def\quest{\refstepcounter{theorem}\paragraph{{\bf Question} \thetheorem}}
\def\definition{\refstepcounter{theorem}\paragraph{{\bf Definition} \thetheorem}}
\newcommand{\bco}{\begin{conjecture}}
\newcommand{\eco}{\end{conjecture}}

\title[${\mathbb P}^1$-gluing for local complete intersections] {``${\mathbb P}^1$-gluing" for local complete intersections}
\author{Mrinal Kanti Das}
\address{Stat-Math Unit, Indian Statistical Institute, 203 B. T. Road, Kolkata 700108 India}
\email{mrinal@isical.ac.in}
\author{Soumi Tikader}
\address{Stat-Math Unit, Indian Statistical Institute, 203 B. T. Road, Kolkata 700108 India}
\email{tsoumi\_r@isical.ac.in,tikadersoumi@gmail.com}
\author{Md. Ali Zinna}
\address{School of Mathematical Sciences, 
National Institute of Science Education and Research Bhubaneswar (HBNI), Odisha 752050, India.
(Present address: Department of Mathematics and Statistics, Indian Institute of Science Education and Research Kolkata, Mohanpur, Nadia - 741246
West Bengal, India)}
\email{zinna2012@gmail.com}

\thanks{}

\date{\today}


\begin{document}

\begin{abstract}
We prove an analogue of the Affine Horrocks' Theorem for local complete intersection ideals
of height $n$ in $R[T]$, where $R$ is a regular domain of dimension $d$, which is essentially of finite type over an infinite perfect field of characteristic unequal to $2$, and $2n\geq d+3$.
\end{abstract}

\maketitle

\begin{center}
{\emph{Dedicated to Professor S. M. Bhatwadekar on his seventieth birthday.}}
\end{center}

\medskip

\hrule

\section{Introduction}
Traditionally, many questions on the theory of algebraic vector bundles on affine varieties (or, more generally, projective modules over commutative Noetherian rings) have been motivated from topology. On the other hand, in nice situations, algebraic vector bundles are deeply connected with local complete intersection ideals of the associated coordinate ring of the variety. Therefore, results on algebraic vector bundles 
inspire one to ask analogous questions on these ideals. For example, the classical 
Quillen-Suslin Theorem \cite{q,su} asserts that any algebraic vector bundle on the affine space ${\mathbb A}^n_k$ 
($k$ field) is trivial. The following conjecture of M. P. Murthy, in some sense,  may be regarded as an analogue of
the Quillen-Suslin Theorem. Recall that the notation $\mu(-)$ stands for the minimal number of generators.

\bco\label{mu1}\cite{mu1}
Let $k$ be a field and let $A=k[T_1,\cdots,T_d]$ be the polynomial ring. Let $n\in {\mathbb N}$ and
$I\subset A$ be an ideal such that $\hh(I)=n=\mu(I/I^2)$. Then $\mu(I)=n$.
\eco

The conjecture is still open in general. The best known answer to this conjecture is the following result due to Mohan Kumar \cite[Theorem 5]{mo} (see also \cite{d2} for a recent result when $k=\overline{\mathbb F}_p$).

\bt\cite{mo}\label{mohan1}
Let $k$ be a field and let $A=k[T_1,\cdots,T_d]$. Let $I\subset A$ be an ideal such that 
$\mu(I/I^2)=n\geq \dd(A/I)+2$. Then $\mu(I)=n$.
\et

Briefly, Mohan Kumar's method goes as follows.
It is well-known that applying a change of variables  one can assume that $I$ 
contains a monic polynomial in $T_d$ with coefficients from $R:=k[T_1,\cdots,T_{d-1}]$. With this 
observation,
Mohan Kumar actually proves the following more general result:

\bt\cite{mo}\label{mohan2}
Let $R$ be a commutative Noetherian ring and $I\subset R[T]$ be an ideal containing a monic polynomial.
Let  $\mu(I/I^2)=n\geq \dd(R[T]/I)+2$. Then there is a projective $R[T]$-module $P$ of rank $n$ and
a surjective $R[T]$-linear map $\alpha:P\sur I$.
\et

Since projective $k[T_1,\cdots,T_d]$-modules are free by the Quillen-Suslin Theorem,  Mohan Kumar
obtains (\ref{mohan1}) as a corollary of (\ref{mohan2}). 
Later, Mandal improves (\ref{mohan2}) in \cite{ma1}, by showing that $P$ can be taken to be free
 (under the same assumptions as in (\ref{mohan2})), and therefore, $\mu(I)=n$. A closer inspection of Mandal's proof 
shows that  he essentially proves the following:

\bt\cite{ma1}
Let $R$ be a commutative Noetherian ring and $I\subset R[T]$ be an ideal containing a monic polynomial.
Let  $I=(f_1,\cdots,f_n)+I^2$, where $n\geq \dd(R[T]/I)+2$. Then there exist $g_1,\cdots,g_n\in I$ such that $I=(g_1,\cdots,g_n)$ with $g_i-f_i\in I^2$. In other words, the $n$ generators of $I/I^2$ can be lifted to a set of $n$ generators of $I$.
\et

We now recall one result on algebraic vector bundles which is closely associated to the Quillen-Suslin Theorem. It is the Affine Horrocks' Theorem \cite{q}: \emph{Let $R$ be any commutative ring and ${\mathcal E}$ be a vector bundle on ${\mathbb A}^1_R$.
If ${\mathcal E}$ extends to a vector bundle on ${\mathbb P}^1_R$, then ${\mathcal E}$ is
extended from $\text{Spec}(R)$.} In particular, (in algebraic terms) we have:

\bt\cite{q}
Let $R$ be a commutative ring and $P$ be a projective $R[T]$-module. Assume that the projective $R(T)$-module $P\otimes_{R[T]} R(T)$ is free. Then $P$ is a free $R[T]$-module. (Here $R(T)$ is the ring obtained from $R[T]$ by inverting all the monic polynomials).
\et

The above theorem may be termed as a \emph{monic inversion principle} 
for (freeness of) projective modules.
From around the time it was proved, the monic inversion principle became a recurrent 
theme in various allied areas. It has been mostly missing from the area of complete intersections.
Motivated by the preceding discussion on Murthy's conjecture, we therefore pose the following natural question.

\medskip

\quest\label{mainq}
\emph{Let $R$ be a commutative Noetherian ring of dimension $d\geq 2$. Let $I\subset R[T]$ be an ideal 
of height $n$ such that $\mu(I/I^2)=n$. Let $I=(f_1,\cdots,f_n)+I^2$ be given. Assume that  
$IR(T)=(G_1,\cdots,G_n)$ with $G_i-f_i\in I^2 R(T)$. Then,
do there exist $F_1,\cdots,F_n\in I$ such that $I=(F_1,\cdots,F_n)$ with $F_i-f_i\in I^2$?}

\medskip

 \rmk
When $d=n=2$, the answer to the above question is negative (see \cite[Example 3.15]{brs1}). However, in this case it follows from \cite[Section 7]{d1} that 
$\mu(I)=2$. When $d=n\geq 3$,  the above question has an affirmative answer if either
$R$ is local, or if $R$ is an affine domain over an algebraically closed field of characteristic zero \cite[Proposition 5.8]{d1}.

\smallskip

We settle Question \ref{mainq} affirmatively in the following form (Theorem \ref{monicinv2} in the text).

\bt
Let $R$ be a regular domain of dimension 
$d$ which is essentially of finite type over an infinite perfect field $k$ of characteristic unequal to $2$. Let $I\subset R[T]$ be an ideal of height $n$ such that $\mu(I/I^2)=n$, where $2n\geq d+3$. 
Let $I=(f_1,\cdots,f_n)+I^2$ be given. Assume that  
$IR(T)=(G_1,\cdots,G_n)$ with $G_i-f_i\in I^2 R(T)$. Then,
 there exist $F_1,\cdots,F_n\in I$ such that $I=(F_1,\cdots,F_n)$ with $F_i-f_i\in I^2$.
\et

As our methods involve quadratic forms and orthogonal groups, we require $2R= R$. 
We now spend a few words on the line of proof. We have $I=(f_1,\cdots,f_n)+I^2$. Let 
$\omega_I:(R[T]/I)^n\sur I/I^2$ be the corresponding surjective map. Applying Nakayama Lemma we can find $h\in I^2$ such that $h-h^2=f_1 g_1+\cdots +f_n g_n$, for some
$g_1,\cdots,g_n\in R[T]$. Now consider the following pointed set 
$$Q'_{2n}(R[T])=\{(x_1,\cdots x_n,y_1,\cdots,y_n,z)\in R[T]^{2n+1}\,|\,\sum_{i=1}^{n} x_iy_i+z^2=1\}$$
with a base point $(0,\cdots,0,0,\cdots,0,1)$. Let $O_{2n+1}(R[T])$ be the orthogonal group preserving the quadratic form $X_1 Y_1+\cdots +X_n Y_n+Z^2$. 
 Let  $EO_{2n+1}(R[T])$ be the elementary subgroup of $O_{2n+1}(R[T])$. We must point out here that
we are considering the elementary subgroup as defined in \cite[Section 3]{vp}. To the pair $(I,\omega_I)$
we associate the $EO_{2n+1}(R[T])$-orbit of $v=(2f_1,\cdots,2f_n,2g_1,\cdots,2g_n,1-2h)$ in
the orbit space $Q'_{2n}(R[T])/EO_{2n+1}(R[T])$. As it turns out, this association does not depend on the choices 
made for $h$ and $g_1,\cdots,g_n$ above. Further, we show that, to lift $\omega_I$ to a surjection
$\varphi: R[T]^n\sur I$, it is enough to prove the following ``quadratic" analogue
of a monic inversion principle of Ravi Rao \cite[Theorem 1.1]{r2}. See Theorem \ref{ravi2}
below.
\remove{Note that this result 
(Theorem \ref{ravi2} in the text) is in a much more general set-up than we need.}

\bt
Let $A$ be a regular domain containing a perfect field $k$ of characteristic $\not = 2$.
Let $n\geq 2$ and $v=(\alpha_1,\cdots,\alpha_n,\beta_1,\cdots,\beta_n,\gamma)\in Q'_{2n}(A[T])$. Suppose that
there is a monic polynomial $f\in A[T]$ such that $v\equiv (0,\cdots,0,1)$ mod $EO_{2n+1}(A[T]_f)$.
Then $v\equiv (0,\cdots,0,1)$ mod $EO_{2n+1}(A[T])$.
\et

\remove{
In this article we use some very basic techniques from the theory of Euler class groups. For
the convenience of the reader we recall these definitions in Section \ref{ECG}, along with some relevant results. Section \ref{bijection} is intermediary in nature where we prepare to leave 
the Euler class group and take a step towards quadratic forms. We do not need any restriction on the
characteristic of the base field in Sections \ref{ECG} and \ref{bijection}.
}

Our methods in this article do not extend to non-smooth set-up (see Example \ref{smb}). On the other hand, lifting the restriction $2n\geq d+3$ (which is concurrent with Mohan Kumar's range 
$n\geq \dd(R[T]/I)+2$) looks to be an enormous task because then it would give a complete solution to Murthy's conjecture, which seems elusive at the moment.

The layout of this article is as follows. In Section \ref{homotopy}, we first recall the definitions of the homotopy orbits of two pointed sets (one of which is $Q'_{2n}(R)$) and record some results which will be used later. Then we  recall the definition of the Euler class group and discuss this group in light of \emph{homotopy} (this treatment is crucial to our paper). Section \ref{bijection} is also a
preparatory section where we prove some technical results which will eventually enable us to translate 
Question \ref{mainq} to a question on quadratic forms.  We do so in Section \ref{Monic} where we prove the main theorem.

\section{Preliminaries}\label{homotopy}

\subsection{Homotopy}
 We first recall the definitions of two pointed sets and their homotopy orbits from \cite{f}.
In this article, by  \emph{`homotopy'} we shall  mean \emph{`naive homotopy'}, as defined below.

\bd
Let $F$ be a functor originating from the category of rings to the category of sets. 
For a given ring $R$, two elements
$F(u_0),F(u_1)\in F(R)$ are said to be homotopic if there is an element $F(u(T))\in F(R[T])$ such that  
$F(u(0))=F(u_0)$ and $F(u(1))=F(u_1)$.
\ed

 Consider the equivalence relation on $F(R)$ \emph{generated} by homotopies
(the relation is easily seen to be reflexive and symmetric but is not transitive in general).
The set of equivalence classes will be denoted by $\pi_0(F(R))$, and an equivalence class will be called as a \emph{homotopy orbit}.

\smallskip

\noindent
\ul{\bf The pointed set $Q'_{2n}(R)$ and its homotopy orbits:}
Let $R$ be any commutative Noetherian ring. Let $n\geq 2$ and
consider the  set 
$$Q'_{2n}(R)=\{(x_1,\cdots x_n,y_1,\cdots,y_n,z)\in R^{2n+1}\,|\,\sum_{i=1}^{n} x_iy_i+z^2=1\}$$
with a base point $(0,\cdots,0,0,\cdots,0,1)$.
Assume that  $2R=R$ and let $O_{2n+1}(R)$ be the
group of orthogonal matrices  preserving the quadratic form 
$\sum_{i=1}^{n} X_i Y_i+Z^2$. Then there is a natural action of $O_{2n+1}(R)$
and its subgroup $SO_{2n+1}(R)$ on the
set $Q'_{2n}(R)$.
Let  $EO_{2n+1}(R)$ be the elementary subgroup of $SO_{2n+1}(R)$ as defined in \cite[Section 3]{vp}, \cite[p. 1503]{va}. 
As $n\geq 2$, the subgroup
$EO_{2n+1}(R)$ is normal in $SO_{2n+1}(R)$ (see \cite[Lemma 4]{vp}). Indeed, the group 
$EO_{2n+1}(R)$ also naturally
acts on the set $Q'_{2n}(R)$. For the convenience of the reader, we recall the definition of 
$EO_{2n+1}(R)$ from \cite{va,vp} below. We explicitly describe the generators of this group by writing out their actions on a vector $(x_1,\cdots,x_n,y_1,\cdots,y_n,z)$. The first three correspond to the
\emph{long root unipotents}, while the last two correspond to the \emph{short root unipotents}, as mentioned in \cite{va,vp}.


\medskip

\definition\label{transvec} {\bf The elementary subgroup $EO_{2n+1}(R)$:} The group  $EO_{2n+1}(R)$ is the subgroup 
of $SO_{2n+1}(R)$ generated by the following \emph{elementary orthogonal transvections:}

 ($1\leq i,j\leq n$, $i\neq j$, and $\lambda\in R$)
\begin{enumerate}
\item
$(x_1,\cdots,x_n,y_1,\cdots,y_n,z)\\
\mapsto (x_1,\cdots,x_{i-1}, x_i+\lambda x_j,x_{i+1},\cdots, y_{j-1},y_j-\lambda y_i,y_{j+1},\cdots, 
y_n,z)$
\item
$(x_1,\cdots,x_n,y_1,\cdots,y_n,z)
\mapsto (x_1,\cdots,x_{i-1},x_i+\lambda y_j,\cdots, x_j-\lambda y_i, x_{j+1},\cdots,y_n,z)$
\item
$(x_1,\cdots,x_n,y_1,\cdots,y_n,z)
\mapsto  (x_1,\cdots, y_{i-1},y_i+\lambda x_j, \cdots,y_j-\lambda x_i,y_{j+1},\cdots, y_n,z)$
\item
$(x_1,\cdots,x_n,y_1,\cdots,y_n,z)
\mapsto (x_1,\cdots,x_{i-1},x_i+2\lambda z-\lambda^2 y_i,x_{i+1},\cdots,y_n,z-\lambda y_i)$
\item
$(x_1,\cdots,x_n,y_1,\cdots,y_n,z)
\mapsto (x_1,\cdots,y_{i-1},y_i+2\lambda z-\lambda^2 x_i,y_{i+1},\cdots,y_n,z-\lambda x_i)$
\end{enumerate}

\medskip

\rmk
Any reader consulting \cite{va,vp} should be cautioned that we are following a different scheme of notations here. The $0\,\textsuperscript{th}$ coordinate of \cite{va,vp} is the last coordinate here.
Also, they use negative indices for the $y$-coordinates. More precisely, their $x_0$ is our $z$, and their $x_{-i}$ is our $y_i$. The reader may also note that the description of $EO_{2n+1}(R)$ by
Calmes-Fasel \cite{cf} as mentioned in \cite[p. 320]{f} is concurrent with the above definition.

\medskip

\rmk\label{vavpet}
It is easy to see from the description of $EO_{2n+1}(R)$ above that
if $\sigma\in EO_{2n+1}(R)$, then $\sigma$ is homotopic to the identity matrix $I_{2n+1}$. 

\medskip

We shall need the following result, which is a consequence of  a much more general result
proved by Stavrova in \cite[Theorem 1.3]{st}. This is an analogue of \cite[Theorem 3.3]{v}.

\bt[Stavrova]\label{stav}
Let $R$ be a regular ring containing a  perfect field $k$ 
with  $\text{Char}(k)\not = 2$. Let $n\geq 2$ and $\tau(T)\in O_{2n+1}(R[T])$ be
such that $\tau(0)=I_{2n+1}$. Then, $\tau(T)\in EO_{2n+1}(R[T])$.
\et

We shall also need the following result proved by Mandal-Mishra \cite[Theorem 4.2]{mm}.

\bt[Mandal-Mishra]\label{manmi}
 Let $R$ be a regular ring containing a  field $k$ 
with  $\text{Char}(k)\not = 2$. Let $H(T)\in Q'_{2n}(R[T])$. Then, 
there exists 
$\tau(T)\in O_{2n+1}(R[T])$ such that $\tau(0)=I_{2n+1}$ and $H(T)\tau(T)=H(0)$.
\et

We are now ready to prove the following theorem. Recall from the beginning of this section that
two elements $(x_1,\cdots x_n,y_1,\cdots,y_n,z)$ and $(x'_1,\cdots x'_n,y'_1,\cdots,y'_n,z')$ 
from $Q'_{2n}(R)$ are homotopic if there is $(f_1(T),\cdots,f_n(T),g_1(T),\cdots,g_n(T),h(T))\in 
Q'_{2n}(R[T])$ such that $f_i(0)=x_i$, \, $f_i(1)=x'_i$, \,\,$g_i(0)=y_i$, \,$g_i(1)=y'_i$ for $i=1,\cdots,n$, and $h(0)=z$ and $h(1)=z'$. Also, $\pi_0(Q'_{2n}(R))$ is the set of homotopy orbits
of $Q'_{2n}(R)$.

\bt\label{uni1}
Let $R$ be a regular ring containing a perfect field $k$ 
with  $\text{Char}(k)\not = 2$. Then,  for any $n\geq 2$
there is a bijection  $$\eta: \pi_0(Q'_{2n}(R))\stackrel{\sim}\lra Q'_{2n}(R)/EO_{2n+1}(R).$$
\et

\proof
Let $u:=(x_1,\cdots x_n,y_1,\cdots,y_n,z)\in Q'_{2n}(R)$. Let $\eta$ be the map which
takes the homotopy orbit of $u$ to its $EO_{2n+1}(R)$-orbit. We have to check first that 
$\eta$ is well-defined. 

Let $u$ and $v$ be 
two representatives of the same homotopy orbit.
By definition, they will be connected by a finite sequence of homotopies. We first consider the simple case, namely, when  $u$ and $v$ are homotopic. In that case,
there exists $H(T)\in Q'_{2n}(R[T])$ such that $H(0)=u$ and $H(1)=v$. By Theorem \ref{manmi}, there
exists $\tau(T)\in O_{2n+1}(R[T])$ such that $\tau(0)=I_{2n+1}$ and $H(T)\tau(T)=H(0)$. Applying
Theorem \ref{stav} we conclude that $\tau(T)\in EO_{2n+1}(R[T])$. Then, 
$\sigma:=\tau(1)\in EO_{2n+1}(R)$ and
we have 
$$v \sigma=H(1)\tau(1)=H(0)=u.$$
Therefore $u$ and $v$ are connected by an element of $EO_{2n+1}(R)$ (namely, $\sigma$),  and they define the same  $EO_{2n+1}(R)$-orbit. To tackle the general case, we can take the product of those $\sigma$ 
obtained for each homotopy.

The map $\eta$ is clearly surjective. To prove the injectivity, let $u,v\in Q'_{2n}(R)$ 
be such that there exists $\sigma\in EO_{2n+1}(R)$ with $v=u \sigma$. By Remark \ref{vavpet},
$\sigma$ is homotopic to $I_{2n+1}$. Therefore, there exists $\Theta(T)\in  EO_{2n+1}(R[T])$ such that
$\Theta(0)=I_{2n+1}$ and $\Theta(1)=\sigma$. Define 
$$H(T):=u\Theta(T)\in Q'_{2n}(R[T]).$$
Then, $H(0)=u$ and $H(1)=v$, showing that $u$ and $v$ are homotopic.
\qed 

 \medskip

The following corollary is now obvious.

\bc\label{uni1c}
Let $R$ be a regular ring containing a perfect field $k$ 
with  $\text{Char}(k)\not = 2$. Then, for $n\geq 2$,
the relation induced by homotopy on $ Q'_{2n}(R)$ is also transitive, and hence an equivalence relation.
\ec

\noindent
\ul{\bf The pointed set $Q_{2n}(R)$ and its homotopy orbits:}
Let $R$ be any commutative Noetherian ring. Let $n\geq 2$ and
consider the  set 
$$Q_{2n}(R)=\{(x_1,\cdots x_n,y_1,\cdots,y_n,z)\in R^{2n+1}\,|\,\sum_{i=1}^{n} x_i y_i=z-z^2\}$$
with a base point $(0,\cdots,0,0,\cdots,0,0)$.
It is proved in \cite{f} that if  $\frac{1}{2}\in R$, then
there is a bijection $\gb_n:Q_{2n}(R)\to Q'_{2n}(R)$ and its inverse $\ga_n:Q'_{2n}(R)\to
Q_{2n}(R)$ given by
\begin{itemize}
\item
$\beta_n(x_1,\cdots x_n,y_1,\cdots,y_n,z)=(2x_1,\cdots 2x_n,2y_1,\cdots,2y_n,1-2z)$
\item
$\ga_n(x_1,\cdots x_n,y_1,\cdots,y_n,z)=\frac{1}{2}(x_1,\cdots x_n,y_1,\cdots,y_n,1-z)$
\end{itemize}

Note that both $\ga_n$ and $\gb_n$ preserve the base points of the respective sets.
They induce bijections between the sets $\pi_0(Q_{2n}(R))$ and $\pi_0(Q'_{2n}(R))$ (will use the same notations).
By transporting the action of $EO_{2n+1}(R)$ on $Q'_{2n}(R)$ through the bijections given above one 
sees that 
$EO_{2n+1}(R)$ also acts on $Q_{2n}(R)$ in the following way:

$$Mv:=\alpha_n((\beta_n(v))M),$$
for $v\in Q_{2n}(R)$ and $M\in EO_{2n+1}(R)$ (with the assumption $\frac{1}{2}\in R$).
Further, note that the bijections $\ga_n$, $\gb_n$ induce bijections between the sets 
$\pi_0(Q_{2n}(R))$ and $\pi_0(Q'_{2n}(R))$. Combining these with Theorem \ref{uni1}, one 
obtains the following result.

\bt\label{uni2}
Let $R$ be a regular ring containing a perfect field $k$ 
with  $\text{Char}(k)\not = 2$. Then,  for any $n\geq 2$
there is a bijection  $$\pi_0(Q_{2n}(R))\stackrel{\sim}\lra Q_{2n}(R)/EO_{2n+1}(R).$$
\et

We shall require the following corollary later.

\bc\label{homoeq}
Let $R$ be a regular ring containing a perfect field $k$ 
with  $\text{Char}(k)\not = 2$. Then, for $n\geq 2$,
the relation induced by homotopy on $ Q_{2n}(R)$ is also transitive, and hence an equivalence relation.
\ec

We shall need a
 ``moving lemma" for $\pi_0(Q_{2n}(R))$ (Lemma \ref{move1} below), where $R$ is any commutative Noetherian ring. This  has also been proved in \cite[Lemma 2.1.5]{af}. However, we prove it here  using the \emph{prime avoidance lemma}, which the reader may find a bit easier to follow. 
We first record  an easy application of the prime avoidance lemma (for a proof, see \cite[Lemma 7.1.4]{ir}), which we shall use a few times. 

\bl\label{EE}
Let $A$ be a commutative Noetherian  ring and $(a_1,\cdots, a_n,a)\in A^{n+1}$. Then there exist 
$\mu_1,\cdots,\mu_n \in A$ such that 
 $\text{ht}(I_{a})\geq n$, where $I=(a_1+a\mu_1,\cdots,a_n+a\mu_n)$. In
other words, if $\mathfrak{p}\in \text{Spec}(A)$ such that  
 $I\subset \mathfrak{p}$ and $a\notin \mathfrak{p}$, then $\text{ht}(\mathfrak{p})\geq n$.
\el

\rmk\label{bertini}
If $A$ is a geometrically reduced affine algebra over an infinite field then Swan's version of Bertini theorem, as given in \cite[Theorem 2.11]{brs2}, implies that $\mu_1,\cdots,\mu_n$ can be so chosen that the ideal $I=(a_1+a\mu_1,\cdots,a_n+a\mu_n)$ has the additional property that $(A/I)_a$ is a geometrically reduced ring.

\bl{(\bf{Moving Lemma})}\label{move1}
Let $R$ be a  commutative Noetherian ring. Let $(a,b,s)=(a_1,\cdots,a_n,b_1,\cdots,b_n,s)\in Q_{2n}(R)$. Then there 
exists $\mu=(\mu_1,\cdots,\mu_n) \in R^n$ such that 
\begin{enumerate}
 \item The row $(a',b',s')=(a_1+\mu_1(1-s)^2,\cdots,a_n+\mu_n(1-s)^2,b_1(1-\mu b^t),\cdots,b_n(1-\mu b^t),
 s+\mu b^t(1-s))\in Q_{2n}(R)$,
 \item $[(a,b,s)]=[(a',b',s')]$ in $\pi_0(Q_{2n}(R))$ and
 \item $\text{ht}(K)\geq n$, where $K=(a_1+\mu_1(1-s)^2,\cdots,a_n+\mu_n(1-s)^2,s+\mu b^t(1-s))$.
\end{enumerate}
\el
\proof
We consider the row $(a_1,\cdots,a_n,(1-s)^2)\in R^{n+1}$. By Lemma \ref{EE} 
there  exist $\mu_1,\cdots,\mu_n \in R$ such that 
 $\text{ht}(I_{(1-s)^2})\geq n$, where $I=(a_1+\mu_1(1-s)^2,\cdots,a_n+\mu_n(1-s)^2)$. In other words, if
 $\mathfrak{p}\in \text{Spec}(R)$ such that $I\subset \mathfrak{p}$ and $(1-s)\notin \mathfrak{p}$, then $\text{ht}(\mathfrak{p})\geq n$.
 
Set $A=a+T(1-s)^2\mu\in R[T]^n$, then an easy computation yields  that 
 $$Ab^t(1- T\mu b^t)=(1-s)(1-T\mu b^t)-(1-s)^2 (1-T\mu b^t)^2.$$
Setting $B=(1-T\mu b^t)b$, it is easy to check that $(A, B, (1-s)(1-T\mu b^t))\in Q_{2n}(R[T])$. 
Then it follows that $(A, B, 1-(1-s)(1-T\mu b^t)) = (A, B, s+T \mu b^t(1-s))\in Q_{2n} (R[T ])$. Thus (1) and (2) are proved. 

Now we have the following relations among the ideals:
\begin{eqnarray*}
I &= & \left(a_1+\mu_1(1-s)^2,\cdots,a_n+\mu_n(1-s)^2\right)\\
& & =\left(a+\mu(1-s)^2, (1-s)(1-\mu b^t)\right) \cap \left(a + \mu(1-s)^2, s+\mu b^t(1-s)\right)\\
& &=\left(a +\mu(1-s)^2 , (1-s)(1-\mu b^t)\right) \cap K
\end{eqnarray*}

Let $\mathfrak{p}\in \text{Spec}(R)$ such that $K\subset \mathfrak{p}$. As $s+\mu b^t(1-s)\in K\subset \mathfrak{p}$, it follows that 
$(1-s)(1-\mu b^t)\notin \mathfrak{p}$ and therefore, $1-s\notin \mathfrak{p}$. Note that $I\subset K\subset \mathfrak{p}$.
Therefore, by the  first paragraph, $\text{ht}(\mathfrak{p})\geq n$. This proves (3).
\qed

\smallskip

\rmk\label{red}
If $R$ is a geometrically reduced affine algebra over an infinite  field then using Swan's Bertini theorem in place of 
Lemma \ref{EE},
one can choose $K$ to have the additional property that either $K=R$ or  
$K$ is  a reduced ideal. 




\subsection{The Euler class groups and homotopy}\label{ECG} There are various definitions of the Euler class groups, depending on the context one is working on. 
Instead of confusing the reader by cluttering this article with each such definition, here we stick to a prototype which is the easiest to understand. We shall then remark on the other cases with suitable references.
 
Let $R$ be a smooth affine domain of dimension $d\geq 2$ over an
infinite perfect field $k$. We recollect definitions of the Euler class
groups from \cite{brs1}. Our emphasis will also be on the definition
of the Euler class
group given by M. V. Nori 
in terms of homotopy (as appeared in \cite{brs1}). In \cite{dk}, the first named author and Manoj Keshari investigated in 
detail the relation between these two equivalent definitions and their consequences. 
We reproduce some of those results in one place for the convenience of the reader.

Let $R$ be as above.
Let $B$ be the set of pairs
$(m,\omega_m)$ where $m$ is a maximal ideal of $R$ and $\omega_m
:(R/m)^d\surj m/m^2$. Let $G$ be the free abelian group generated by
$B$.  Let $J=m_{1}\cap \cdots \cap m_r$, where $m_i$ are distinct maximal
ideals of $R$. Any $\omega_{J}:(R/J)^d\surj J/J^2$ induces surjections
$\omega_i :(R/m_i)^d\surj m_i/m_i^2$ for each $i$. We associate
$(J,\omega_J):= \sum_{1}^{r}(m_i,\omega_i)\in G$.  Now, we have the following definitions:

\medskip

\definition\label{nori} (Nori) Let $S$ be the set of elements
$(I(1),\omega(1))-(I(0),\omega(0))$ of $G$ where (i) $I\subset R[T]$
is a local complete intersection ideal of height $d$; (ii) Both $I(0)$
and $I(1)$ are reduced ideals of height $d$; (iii) $\omega(0)$ and
$\omega(1)$ are induced by $\omega: (R[T]/I)^d\surj I/I^2$. Let $H$ be
the subgroup generated by $S$. The $d\,\textsuperscript{th}$ Euler class group $E^d(R)$ is defined
as $E^d(R):=G/H$.

\medskip

\begin{definition} \label{brs}(Bhatwadekar-Sridharan) Let $S_1$ be
the set of elements $(J,\omega_J)$ of $G$ for
which $\omega_J$ has a lift to a surjection $\theta: R^d\surj J$ and $H_1$ be the 
subgroup of $G$ generated by $S_1$ .  The
Euler class group $E^d(R)$ is defined as $E^d(R):=G/H_1$.
\end{definition}

\medskip

\rmk
We shall refer to the elements of the Euler class group as \emph{Euler cycles}.

\medskip

\rmk
The above definitions appear to be  slightly different than the ones given in \cite{brs1}. 
However, note that if $(J,\omega_J)\in S$ (resp. $S_1$) and if $\ol\sigma\in E_d(R/J)$, then
the element $(J,\omega_J\ol\sigma)$ is also in  $S$ (resp. $S_1$). For details, see \cite[Remark 5.4]{dk}
and \cite[Proposition 2.2]{dz}. 

\medskip

\rmk
Bhatwadekar-Sridharan proved (see \cite[Remark 4.6]{brs1}) that $H=H_1$ and therefore the above definitions of the Euler class group are equivalent.

\medskip




\smallskip

The following theorem collects a few results in one place
(see \cite[4.11]{brs1}, \cite[4.2]{k}, \cite[Theorem 5.13]{dk} for details).

\bt\label{zero}
Let $R$ be a smooth affine domain of dimension $d\geq 2$ over an
infinite perfect field $k$. Let $J\subset R$ be a reduced ideal of height $d$ and
$\omega_{J}:(R/J)^d\surj J/J^2$ be a surjection. Then, the following are equivalent:
\begin{enumerate}
\item
The image of $(J,\omega_J)=0$ in $E^d(R)$
\item
$\omega_J$ can be lifted to a surjection $\theta: R^d\sur J$.
\item
$(J,\omega_J)=(I(0),\omega(0))-(I(1),\omega(1))$ 
in $G$
where  (i) $I\subset R[T]$ is a local complete intersection ideal of height $d$; (ii) Both $I(0)$ and $I(1)$ are
reduced ideals of height $d$, and (iii) $\omega(0)$ and $\omega(1)$ are induced by $\omega: (R[T]/I)^d\sur I/I^2$.  
\end{enumerate}
\et

A series of remarks are in order.

\remove{WRONG REMARK----------
\rmk\label{crucial}
Let $(J,\omega_J), (J',\omega_{J'})$ be such that both $J,J'$ are reduced ideals and they represent the same element in $E^d(R)$. In other words, $(J,\omega_J)- (J',\omega_{J'})\in H_1$
(where $H_1$ is as in (\ref{brs}). It is easy to see from the proof of \cite[4.1]{k} that 
 $(J,\omega_J)- (J',\omega_{J'})$ actually
belongs to $S_1$ (where $S_1$ is as in (\ref{brs})). This observation will be crucially used in the proof of Proposition \ref{map1}.}

\smallskip

\rmk
Let $J\subset R$ be an ideal of height $d$ which is not necessarily reduced and let 
$\omega_{J}:(R/J)^d\surj J/J^2$ be a surjection. Then also one can associate an element 
$(J,\omega_J)$ in $E^d(R)$ and prove the above theorem for $(J,\omega_J)$. See \cite[Remark 4.16]{brs1} for details.

\smallskip

\rmk\label{ess}
All the definitions and results in this subsection can be easily extended to the case when  $R$ is a regular domain of dimension $d\geq 2$ which is essentially of finite type over  an infinite perfect 
field $k$. The key result on which this entire subsection depends (including the equivalence of the two definitions of $E^d(R)$) is the \emph{``homotopy theorem"} of 
Bhatwadekar-Sridharan \cite[Theorem 3.8]{brs1}. This homotopy theorem, which was proved for smooth affine domains over an infinite perfect field, has an obvious extension to regular domains essentially of finite type over such fields (for instance, see \cite[Theorem 4.13]{bk}).

\smallskip

\rmk\label{lower}
All the definitions and results in this subsection can also be extended  to
a much more relaxed range. Let $R$ be a regular domain of dimension $d$, which is essentially of finite type over  an infinite perfect field $k$.
Let $n$ be an integer such that $2n\geq d+3$. Then 
the $n\textsuperscript{th}$ Euler class group $E^n(R)$  has been defined in \cite{brs4}.  
To extend the results of this section one has to use \cite[Theorem 4.13]{bk}, which is
a generalization  of the homotopy theorem of Bhatwadekar-Sridharan mentioned above.
We do not work 
out the details here as  the process is routine. 

\section{A bijection}\label{bijection}
In this section we need to prove some more preparatory results in order to be able to prove the main theorem of this article in the next section. 
 
As we did with the Euler class groups in the previous section, here also we shall work out the details
when $R$ is a smooth affine domain of dimension $d\geq 2$  over  an infinite perfect 
field $k$. We shall then remark on the extensions to the other cases.

\smallskip

\rmk
Throughout this section we do not need any assumption on the characteristic of the base field $k$.
This section is all about the relation between the Euler class group $E^d(R)$ and $\pi_0(Q_{2d}(R))$.
Note that, in order to talk about $\pi_0(Q_{2d}(R))$ we do not need $2$ to be invertible.

\medskip

\bd\label{deftheta}
\ul{\bf A set-theoretic map:}
We first define a set-theoretic map from the Euler class group $E^d(R)$ to $\pi_0(Q_{2d}(R))$.
By \cite[Remark 4.14]{brs1} we know that an arbitrary element of $E^d(R)$ can be represented by a single Euler cycle 
$(J,\omega_J)$, where $J$ is a reduced ideal of height $d$. Now $\omega_J:(R/J)^d\sur J/J^2$ is given by $J=(a_1,\cdots,a_d)+J^2$, for
some $a_1,\cdots,a_d\in J$.
Applying the  Nakayama Lemma one obtains $s\in J^2$
such that  $J=(a_1,\cdots,a_d,s)$ 
 with $s-s^2=a_1b_1+\cdots +a_d b_d$  for some $b_1,\cdots,b_d\in R$ (see \cite{mo1} for a proof).
We  associate to $(J,\omega_J)$ the homotopy orbit  $[(a_1,\cdots,a_d,b_1,\cdots,b_d,s)]$
in $\pi_0(Q_{2d}(R))$. 
\ed

\medskip

\rmk
For an ideal $I$ with a given set of generators of $I/I^2$, the idea of associating the whole data
(as above) to an affine quadric,  perhaps goes back to Mohan Kumar and Nori (to the best of our knowledge). We refer to
\cite[Section 17]{sw} for an  exposition of their method. The novelty of the approach in \cite{f} is to associate the information up to naive homotopy, which we borrow in this article.

\medskip

\rmk
In the proof of Proposition \ref{map1} below, we shall use \cite[Theorem 2.0.2]{f}. We must point out to the reader that the gap in the said paper of Fasel was precisely in \cite[Lemma 3.2.3]{f} (see 
\cite{f2}), whereas the proofs in Section 2 of \emph{Op. Cit.} are sound. 

\medskip


\bp\label{map1}
Let $R$ be a smooth affine domain of dimension $d\geq 2$  over  an infinite perfect 
field $k$.
The association $(J,\omega_J)\mapsto [(a_1,\cdots,a_d,b_1,\cdots,b_d,s)]$ is well defined and gives rise to a set-theoretic map $\theta_d: E^d(R)\to \pi_0(Q_{2d}(R))$. The map $\theta_d$ takes the
trivial Euler cycle to the homotopy orbit of the base point $(0,\cdots,0)$ of $Q_{2d}(R)$.
\ep

\proof
We need to check the following:

\begin{enumerate}
\item
If $\omega_J$ is also given by $J=(\ga_1,\cdots,\ga_d)+J^2$ and if $\tau\in J^2$ is such that 
$\tau-\tau^2=\ga_1\gb_1+\cdots +\ga_d\gb_d$, then $[(a_1,\cdots,a_d,b_1,\cdots,b_d,s)]=
[(\ga_1,\cdots,\ga_d,\gb_1,\cdots,\gb_d,\tau)]$
in $\pi_0(Q_{2d}(R))$.
\item
If $\ol\sigma\in E_d(R/J)$, then
the image of $(J,\omega_J\ol\sigma)$ in $Q_{2d}(R)$ is homotopic to 
the image of $(J,\omega_J)$.
\item
If $(J,\omega_J)=(J',\omega_{J'})$  in $E^d(R)$ (where both $J, J'$ are reduced ideals
of height $d$), then their images
are homotopic in $Q_{2d}(R)$. 
\end{enumerate}

\noindent
\emph{Proof of (1)} : This has been proved  in \cite[Theorem 2.0.2]{f}. 

\noindent
\emph{Proof of (2)} :
Suppose that $(J,\omega_J)$ is given by $J=(a_1,\cdots,a_d)+J^2$, and 
$s\in J^2$
be such that  $J=(a_1,\cdots,a_d,s)$ 
 with $s-s^2=a_1b_1+\cdots +a_d b_d$  for some $b_1,\cdots,b_d\in R$.
Let $\sigma\in E_d(R)$ be a lift of $\ol\sigma$ and write
$(a_1,\cdots,a_d)\sigma=(\ga_1,\cdots,\ga_d)$. Then $J=(\ga_1,\cdots,\ga_d)+J^2$
and $J=(\ga_1,\cdots,\ga_d,s)$.
If we write $(b_1,\cdots,b_d)(\sigma^{-1})^t=(\gb_1,\cdots,\gb_d)$ (here $t$ stands for transpose), then it is easy to
see that $s(1-s)=\ga_1\gb_1+\cdots +\ga_d\gb_d$. Now, note that

$$\lambda=
\begin{pmatrix}
  
  \sigma & 0 & 0 \\
  0 & (\sigma^{-1})^t & 0\\
   0 & 0 & 1
\end{pmatrix}  \in E_{2d+1}(R)$$

and $(a_1,\cdots,a_d,b_1,\cdots,b_d,s)\lambda=(\ga_1,\cdots,\ga_d,\gb_1,\cdots,\gb_d,s)$. Since
elementary matrices are homotopic to identity, we are done in this case.

\medskip

\noindent
\emph{Proof of (3)} : We break this proof into two steps.

\noindent
\emph{Step 1.}
We have
$(J,\omega_J)=(J',\omega_{J'})\in E^d(R)=G/H_1$, where $H_1$ is as in Definition \ref{brs}. 
It then follows that
$$(J,\omega_J)+\sum_{i=1}^{r}(K_i,\omega_{K_i})=(J',\omega_{J'})+\sum_{j=r+1}^{s}(K'_j,\omega_{K'_j})$$
in $G$, where all the  $(K_i,\omega_{K_i})$ and $(K'_j,\omega_{K'_j})$ are in $S_1$
(where $S_1$ is as in Definition \ref{brs}). 

Adapting the proof of \cite[4.11]{brs1}, if necessary,  we can change the above equation to obtain a new one where the ideals appearing on the left are mutually comaximal (consequently, so are the ideals on the right). Therefore, without loss of generality, we may assume that 
$J,K_1,\cdots,K_r$ are mutually comaximal (and so are $J',K'_{r+1},\cdots,K'_s$) and we 
have $J\cap K_1\cap \cdots\cap K_r=J'\cap K'_{r+1}\cap \cdots\cap K'_s$. Let us write $K=K_1\cap \cdots\cap K_r$ and $K'=K'_{r+1}\cap \cdots\cap K'_s$. Also, let $\omega_K$ be the surjection
$(R/K)^d\sur K/K^2$ induced by $\omega_{K_1},\cdots,\omega_{K_r}$. Similarly, let $\omega_{K'}$ be the surjection
$(R/K')^d\sur K'/K'^2$ induced by $\omega_{K'_{r+1}},\cdots,\omega_{K'_s}$. Summing up, we have:
\begin{enumerate}
\item
$(J,\omega_J)+(K,\omega_K)=(J',\omega_{J'})+(K',\omega_{K'})$ in $G$;
\item
$J+K=R=J'+K'$;
\item
$J\cap K=J'\cap K'$.
\end{enumerate}

\smallskip

\noindent
\emph{Step 2.}
Assume that $\omega_J$ is induced by $J=(a_1,\cdots,a_d)+J^2$ and let $\omega_K$ be induced by 
$K=(c_1,\cdots,c_d)$. They will together induce $\omega_{J\cap K}:(R/J\cap K)^d\sur (J\cap K)/(J\cap K)^2$,
say, given by $J\cap K=(\gb_1,\cdots,\gb_d)+(J\cap K)^2$.

Because of (2) above, we are now free to apply elementary transformations.
Applying elementary transformations on $(c_1,\cdots,c_d)$, if necessary, we may assume by 
\cite[Lemma 3]{rs1} that $\hh(c_1,\cdots,c_{d-1})=d-1$ and 
$J+(c_1,\cdots,c_{d-1})=R$ (Note that if we apply $\sigma\in E_d(R)$ on 
$(c_1,\cdots,c_d)$, we have to apply $\sigma$ on $(a_1,\cdots,a_d)$ as well to retain the
relations  and equations). Consider the ideal $L=(c_1,\cdots,c_{d-1},(1-c_d)T+c_d)$ in $R[T]$.
Write $I=L\cap J[T]$. Using the Chinese Remainder Theorem we can then find 
$f_1,\cdots,f_d\in I$ such that: 
\begin{enumerate}
\item[(a)]
$I=(f_1,\cdots,f_d)+I^2$.
\item[(b)]
$f_i=c_i$ mod $L^2$ for $i=1,\cdots,d-1$ and $f_d=(1-c_d)T+c_d$
mod $L^2$. 
\item[(c)]
 $f_i=a_i$ mod $J[T]^2$ for $i=1,\cdots,d$.
\end{enumerate}

Let $\omega:(R[T]/I)^d\sur I/I^2$ be the surjection corresponding to $f_1,\cdots,f_d$. We then have,
$I(0)=J\cap K$, $I(1)=J$. From (b) we get $f_i(0)=c_i$ mod $K^2$. On the other hand, from (c) we get $f_i(0)=a_i$ mod $J^2$. Combining, we observe that $f_i(0)\equiv \gb_i$ mod $(J\cap K)^2$. In other words, $\omega_{J\cap K}$ is the same as  $\omega(0)$. Also, from (c), we obtain $f_i(1)\equiv a_i$ mod $J^2$,
implying that $\omega_J$ is the same as $\omega(1)$. Therefore, by  Lemma \ref{homlem} proved below,
the images of $(J,\omega_J)$ and $(J\cap K,\omega_{J\cap K})$ are the same in $\pi_0(Q_{2d}(R))$.

Following exactly the same procedure, as above, we can see that the images of $(J',\omega_{J'})$
and $(J'\cap K',\omega_{J'\cap K'})$ are the same in $\pi_0(Q_{2d}(R))$.  Since 
$J\cap K=J'\cap K'$, and $\omega_{J\cap K}=\omega_{J'\cap K'}$, it follows that $(J,\omega_J)$ and $(J',\omega_{J'})$ have the same image in $\pi_0(Q_{2d}(R))$.
This completes the proof of the proposition.
\qed

\bl\label{homlem}
Let $R$ be a smooth affine domain of dimension $d\geq 2$ over an infinite perfect field $k$.
Let  $I\subset R[T]$ be an ideal of height $d$ such that both $I(0)$ and $I(1)$ are ideals of height $d$ in $R$. Assume that there is a surjection $\omega:(R[T]/I)^d\sur I/I^2$.  Then, the images of
$(I(0),\omega(0))$ and $(I(1),\omega(1))$ in $\pi_0(Q_{2d}(R))$ are the same.
\el

\proof
Let $\omega:(R[T]/I)^d\sur I/I^2$ be given by $I=(f_1(T),\cdots,f_d(T))+I^2$. Then, $\omega(0)$ is given by $I(0)=(f_1(0),\cdots,f_d(0))+I(0)^2$ and $\omega(1)$ is given by $I(1)=(f_1(1),\cdots,f_d(1))+I(1)^2$.

There exist $h(T)\in I^2$ and $g_1(T),\cdots, g_d(T)\in R[T]$ such that:
\begin{enumerate}
\item
$I=(f_1(T),\cdots,f_d(T),h(T))$;
\item
$h(T)-h(T)^2=f_1(T)g_1(T)+\cdots +f_d(T)g_d(T)$. 
\end{enumerate}
Then the $2d+1$-tuple
$(f_1(T),\cdots,f_d(T),g_1(T),\cdots, g_d(T),h(T))\in Q_{2d}(R[T])$.

Now  $h(0)\in I(0)^2$ with $h(0)-h(0)^2=f_1(0)g_1(0)+\cdots+f_d(0)g_d(0)$.
 Similarly, we have $h(1)\in I(1)^2$ with $h(1)-h(1)^2=f_1(1)g_1(1)+\cdots+f_d(1)g_d(1)$. Therefore it is easy to see
 that $(f_1,\cdots,f_d,g_1,\cdots,g_d,h)\in Q_{2d}(R[T ])$ is the required homotopy for the images of
$(I(0),\omega(0))$ and $(I(1),\omega(1))$ in $Q_{2d}(R)$.
 This concludes the proof.
\qed

\medskip

\rmk
So far we have worked with Euler cycles represented by reduced ideals. Now let $J$ be an ideal which is not reduced and $\omega_J:(R/J)^d\sur J/J^2$ be a surjection. As indicated in \cite[4.16]{brs1}, using Swan's Bertini theorem we can find a reduced ideal $K$ of height $d$ and elements $a_1,\cdots,a_d$ such that: (i) $J\cap K=(a_1,\cdots,a_d)$; (ii) $J+K=R$; (iii) the images of $a_1,\cdots,a_d$ induce $\omega_J$ (for a proof, see \cite[Lemma 2.7, Remark 2.8]{drs1}). Let 
$\omega_K:(R/K)^d\sur K/K^2$ be the surjection induced by $a_1,\cdots,a_d$. We may apply the same procedure again and find a reduced ideal $L$ of height $d$ such that: (iv) $K\cap L=(b_1,\cdots,b_d)$;
(v) $L+K\cap J=R$; (vi) 
$b_i-a_i\in K^2$ for $i=1,\cdots,d$. Let $\omega_L:(R/L)^d\sur L/L^2$ be the surjection induced by 
$b_1,\cdots,b_d$. One then associates $(J,\omega_J):=(L,\omega_L)$ in $E^d(R)$. It can be easily checked that this association is well-defined and $(J,\omega_J)$ also satisfies the calculus of Euler cycles represented by reduced ideals. Now, to the data $J=(a_1,\cdots,a_d)+J^2$ we can associate an
element of $\pi_0(Q_{2d}(R))$ (exactly as we did in Definition \ref{deftheta}). On the other hand, from $L=(b_1,\cdots,b_d)+L^2$, we shall obtain $\theta_d((L,\omega_L))\in \pi_0(Q_{2d}(R))$.
We now prove:

\bp
With notations as above, the element of $\pi_0(Q_{2d}(R))$ associated to $J=(a_1,\cdots,a_d)+J^2$ is the same as $\theta_d((L,\omega_L))$. 
\ep

\proof
We first note that statements (1) and (2) of Proposition \ref{map1}  are also true for the pair $(J,\omega_J)$.

The idea of proof of this proposition is essentially contained in the proof of 
Proposition \ref{map1} (3) (\emph{Step 2}). Therefore, we shall only give a sketch. As $K\cap L=(b_1,\cdots,b_d)$, we can easily construct an ideal $I\subset R[T]$ and a surjection 
$\omega:(R[T]/I)^d\sur I/I^2$ such that $I(0)=J$, $I(1)=J\cap K\cap L$, $\omega(0)=\omega_J$,
and $\omega(1)=\omega_{J\cap K\cap L}$. On the other hand, as $J\cap K=(a_1,\cdots,a_d)$,
we can construct an ideal $I'\subset R[T]$ and a surjection 
$\omega':(R[T]/I)^d\sur I'/{I'}^2$ such that $I'(0)=L$, $I'(1)=J\cap K\cap L$, $\omega'(0)=\omega_L$,
and $\omega'(1)=\omega_{J\cap K\cap L}$. We can now apply Lemma \ref{homlem} to conclude the proof.
\qed


\medskip

\rmk
Proposition \ref{map1} has also been proved in \cite[Proposition 3.1.9]{af}, 
\cite[Section 6.1]{mm} with the additional assumption that 
$\text{Char}(k)\neq 2$. We do not need that assumption and our line of proof is different. 

\medskip

We now  prove that the set-theoretic map $\theta_d: E^d(R)\to \pi_0(Q_{2d}(R))$
is a bijection. 

\bt\label{biject}
Let $R$ be a smooth affine domain of dimension $d\geq 2$  over  an infinite perfect 
field $k$. 
The set-theoretic map $\theta_d: E^d(R)\to \pi_0(Q_{2d}(R))$ is a bijection.
\et

\proof
Let $v=(a_1,\cdots,a_d,b_1,\cdots,b_d,s)\in Q_{2d}(R)$. Then the ideal $I(v):=(a_1,\cdots,a_d,s)$
of $R$ need not be of height $d$. However, we may apply Lemma \ref{move1} 
 to obtain $v'=(a'_1,\cdots,a'_d,b'_1,\cdots,b'_d,s')$ in the 
same homotopy class of $v$ such that the ideal  $K=(a'_1,\cdots,a'_d,s')$ has height $\geq d$. 
Assume  that $K$ is proper.
We have $K=(a'_1,\cdots,a'_d)+K^2$. If $\omega_K:(R/K)^d\sur K/K^2$ is the corresponding map,
then it follows that the image of $(K,\omega_K)$ under $\theta_d$ is $[v']=[v]$ in  $\pi_0(Q_{2d}(R))$.
On the other hand, if $K=R$, then the row $(a'_1,\cdots,a'_d,s')$ is unimodular and therefore 
there exist $\ga_1,\cdots,\ga_d,\gb\in R$ such that 
$$\ga_1 a'_1+\cdots + \ga_d a'_d+\gb s'=1, \text{ and therefore,}$$
$$(1-s')(\ga_1 a'_1+\cdots + \ga_d a'_d)+\gb(s'-{s'}^2)=1-s'.$$
As $s'-{s'}^2$ is in the ideal $(a'_1,\cdots,a'_d)$, it follows that there exist $\gl_1,\cdots,\gl_d\in R$ such that $\gl_1 a'_1+\cdots +\gl_d a'_d=1-s'$. We can apply elementary orthogonal transformation of type 5 in Definition \ref{transvec} and change $(a'_1,\cdots,a'_d,b'_1,\cdots,b'_d,s')$ to
$(a'_1,\cdots,a'_d,b''_1,\cdots,b''_d,1)$. The latter is clearly homotopic to 
$(0,\cdots,0,0,\cdots,0,1)$. By \cite[Lemma 5.3]{mm}, the orbits $[(0,\cdots,0,0,\cdots,0,1)]$ and
$[(0,\cdots,0,0,\cdots,0,0)]$ are the same in $\pi_0(Q_{2d}(R))$. This trivial orbit has preimage in
$E^d(R)$.
Therefore, $\theta_d$ is surjective.

The rest of the proof is devoted to  proving that $\theta_d$ is injective. Let $(J,\omega_J)$ and $(J',\omega_{J'})$ be elements of $E^d(R)$ be such that
$\theta_d((J,\omega_J))=\theta_d((J',\omega_{J'}))$.
Let $\omega_J$ be given by $J=(a_1,\cdots,a_d)+J^2$. As $\hh(J)=d$, applying Lemma \ref{EE} if necessary, we may assume that $\hh(a_1,\cdots,a_d)=d$. Now there exists $s\in J^2$  such
that $J=(a_1,\cdots,a_d,s)$ with $s-s^2=a_1b_1+\cdots+a_db_d$ for some $b_1,\cdots,b_d\in R$. Similarly,
$\omega_{J'}$ is given by $J'=(a_1',\cdots,a_d')+J'^2$ with $\hh(a_1',\cdots,a_d')=d$. There exists 
$s'\in J'^2$ be such that $J'=(a_1',\cdots,a_d',s')$  with
$s'-s'^2=a_1'b_1'+\cdots+a_d'b_d'$ for some $b_1',\cdots,b_d'\in R$.


We now assume that  
$$\theta_d ((J,\omega_J))=[(a_1,\cdots,a_d,b_1,\cdots, b_d,s)]=[(a'_1,\cdots,a'_d,b'_1,\cdots, b'_d,s')]
=\theta_d ((J',\omega_{J'}))$$
in $\pi_0(Q_{2d}(R))$. Applying Corollary \ref{homoeq} we have $V=(f_1,\cdots,f_d,g_1,\cdots, g_d,h)\in Q_{2d}(R[T])$
such that $V(0)=(a_1,\cdots,a_d,b_1,\cdots, b_d,s)$ and $V(1)=(a'_1,\cdots,a'_d,b'_1,\cdots, b'_d,s')$.
If we consider the ideal $I=(f_1,\cdots,f_d,h)$ of $R[T]$ then we have $I=(f_1,\cdots,f_d)+I^2$. Let
$\omega_I:(R[T]/I)^d\sur I/I^2$ denote the corresponding surjection.
However, the height of $I$ need not be  $d$, although both $I(0)\,\,(=J)$ and $I(1)\,\,(=J')$ have height $d$.  

As both $\hh((a_1,\cdots,a_d)=d=\hh(a_1',\cdots,a_d')$, it follows that 
$$\hh(f_1,\cdots,f_d,T(T-1))=d+1 \,\,\,\,\,\,\,\,\,\,\,\,\,\,\,(\ast)$$ 
Consider $(f_1,\cdots,f_d,(T^2-T)h^2)\in R[T]^{d+1}$. By Lemma \ref{EE},
there exist $\mu_1,\cdots,\mu_d\in R[T]$ such that 
$\hh((F_1,\cdots,F_d)_{h^2(T^2-T)})\geq d,$ where
 $F_i=f_i+\mu_ih^2(T^2-T)$, for $i=1,\cdots,d$. 
Note that we have
$I=(F_1,\cdots,F_d)+(h),$ and $(h)\subset I^2$. Applying \cite[2.11]{brs3}, there exists $e\in (h)$ such that
$I=(F_1,\cdots,F_d,e)$ where $e-e^2\in (F_1,\cdots,F_d)$. 
We now take $K=(F_1,\cdots,F_d,1-e)$ and write $\omega_K: (R[T]/K)^d\sur K/K^2$ for the
corresponding surjection. We record that $I\cap K=(F_1,\cdots,F_d)$ in $R[T]$.

Let $P\in \text{Spec}(R[T])$ be such that $K\subseteq P$. Then, as $e\in (h)$ and $1-e\in K$, we see that $h\not\in P$. If $T^2-T\not\in P$, 
then $\hh(P)\geq d$.
If $T^2-T\in P$, then by ($\ast$) above,  $\hh(P)\geq d+1$. In any case, $\hh(K)\geq d$. 
Note that $$K(0)\cap I(0)=K(0)\cap J=(F_1(0),\cdots,F_d(0))=(a_1,\cdots,a_d),$$
$$K(1)\cap I(1)=K(1)\cap J'=(F_1(1),\cdots,F_d(1))=(a_1',\cdots,a_d').$$

As the height of each of  the ideals involved here is $d$, we have
$$(J,\omega_J)+(K(0),\omega_{K(0)})=0=(J',\omega_{J'})+(K(1),\omega_{K(1)}) \text{ in } E^d(R),$$
where $\omega_{K(0)}$ is induced by $a_1,\cdots,a_d$, and $\omega_{K(1)}$ is induced by 
$a_1',\cdots,a_d'$.

 Therefore, 
$(J,\omega_J)-(J',\omega_{J'})=(K(1),\omega_{K(1)})-(K(0),\omega_{K(0)})\in H$
(where $H$ is as in (\ref{nori})) and consequently, $(J,\omega_J)=(J',\omega_{J'})$ in $E^d(R)$.
This completes the proof.
\qed


\medskip

\rmk
In Proposition  \ref{map1} and Theorem \ref{biject} we can take $R$ to be a regular domain of dimension $d$ which is essentially of finite type over  an infinite  perfect field $k$. Moreover,
following the same line of arguments as in Proposition \ref{map1} and Theorem \ref{biject}, with suitable modifications, any diligent reader can see that there is a similar set-theoretic map
$\theta_n:E^n(R)\to \pi_0(Q_{2n}(R))$ (where $R$ is a regular domain of dimension $d$ which is essentially of finite type over  an infinite  perfect field $k$, and $n$ is an integer such that 
$2n\geq d+3$), and prove the following theorem.

\bt\label{biject2}
Let $R$ be a regular domain of dimension $d$ which is essentially of finite type over  an infinite  perfect field $k$. Let $n$ be an integer such that $2n\geq d+3$.
Then there is a bijection $\theta_n:E^n(R)\lra \pi_0(Q_{2n}(R))$.
\et

\medskip

\rmk
We need to clarify certain points here, as follows.
\begin{enumerate}
\item
Theorem \ref{biject2} has been proved independently in \cite[Theorem 3.1.13]{af} (but with quite a few additional assumptions),  and in \cite[Theorem 6.15]{mm} with the restriction that 
$\text{Char}(k)\neq 2$. Also, one can easily see that our approach and arguments are entirely different from these two articles.
\remove{
\item
Unlike \cite{af,mm}, we do not make any assumption on  the characteristic of $k$.
Note that, to talk \emph{only} about $\pi_0(Q_{2n}(R))$, one does not need $2$ to be invertible in $R$.}
\item
The case $d=n=2$ has not been covered in \cite{af,mm}, which we do in Theorem \ref{biject}. This case is crucial to our paper \cite{dtz}. Let $R$ be a smooth affine algebra of dimension $d\geq 2$ over 
$\br$. In \cite{dtz} we use Theorem \ref{biject} to define a morphism 
$E^d(R)\lra Um_{d+1}(R)/E_{d+1}(R)$ 
which enables us to prove a series of results. For details we refer to \cite{dtz}. 
\item
Because of (2), we decided to prove Theorem \ref{biject} in its present form.
\item
The bijection $\theta_n:E^n(R)\lra \pi_0(Q_{2n}(R))$ from Theorem \ref{biject2} induces a group 
structure on $\pi_0(Q_{2n}(R))$.
\end{enumerate}

\rmk
Let $R$ be a commutative Noetherian ring of dimension $d\geq 2$ which is not necessarily regular. The $d$-th Euler class group $E^d(R)$ has been defined in this case as well 
(see \cite{brs3}). Further,  following the same set of arguments as above, one can define a well-defined set-theoretic map $\theta:E^d(R)\longrightarrow \pi_0(Q_{2d}(R))$. This map is surjective but is not injective in general. In fact, even
if $R$ is an affine algebra over an algebraically closed field, $\theta$ may not be injective. The issue here is smoothness. We elaborate with an example.

\medskip

The following example of Bhatwadekar, based on an example constructed by Bhatwadekar, Mohan Kumar and Srinivas \cite[Example 6.4]{brs1}, appeared in \cite{dk}. We recall it verbatim.

\medskip

\begin{example}\label{smb}(Bhatwadekar)
 Let
$$B=\frac{\BC [X,Y,Z,W]}{(X^5+Y^5+Z^5+W^5)},$$
as in \cite[Example 6.4]{brs1}. 
Then $B$ is a graded
normal affine domain over $\BC$ of dimension $3$, having an isolated
singularity at the origin. Let $F(B)$ be the subgroup of $\wt{K}_0(B)$
generated by all elements of the type $[P]-[P^\ast]$, where $P$ is a
finitely generated projective $B$-module. As $B$ is graded,
$\text{Pic}(B)=0$. Therefore, by \cite[6.1]{brs1}
$F(B)=F^3K_0(B)$. Since $\text{Proj}(B)$ is a smooth surface of degree
$5$ in $\BP^3$, it follows from a result of Srinivas that $F(B)=
F^3K_0(B)\not = 0$. Therefore, there exists a projective $B$-module
$P$ of rank $3$ with trivial determinant such that $[P] - [P^{\ast}]$
is a nonzero element of $F(B)$.  This implies that $P$ does not have a
unimodular element.  We now consider the ring homomorphism $f: B
\rightarrow B[T]$ given by $f(x) = xT, f(y) = yT, f(z) = zT, f(w) =
wT$.  We regard $B[T]$ as a $B$-module through this map and $Q = P
\otimes_B B[T]$.  Then it is easy to see that $Q/TQ$ is free and
$Q/(T-1)Q = P$.  Therefore, $Q$ is a projective $B[T]$-module which is
not extended from $B$. Now consider a surjection $\alpha : Q\surj I$
where $I\subset B[T]$ is an ideal of height $3$. Fix an isomorphism
$\chi:B[T]\simeq \wedge^3(Q)$. Note that $Q/IQ$ is a free
$B[T]/I$-module. We choose an isomorphism $\sigma:(B[T]/I)^3\simeq
Q/IQ$ such that $\wedge^3\sigma= \chi\otimes B[T]/I$. Composing
$\sigma$ and $\alpha\otimes B[T]/I$, we obtain a surjection
$\omega:(B[T]/I)^3\surj I/I^2$.  It is now easy to see that
$(I(0),\omega(0))=0$ in $E^3(B)$ (as $Q/TQ$ is free), whereas
$(I(1),\omega(1))=e(Q/(T-1)Q,\chi(1))=e(P,\chi(1))$ in $E^3(B)$ cannot
be trivial (as $P$ does not have a unimodular element).
$\qed$ 
\end{example}

\medskip

\rmk
The above example shows that there are distinct elements in $E^3(B)$ whose images in 
$\pi_0(Q_{6}(B))$ are the same.

\section{Main theorem: The  monic inversion principle}\label{Monic}
In this final section we address the following question raised in the introduction.
Recall that $R(T)$ is the ring obtained from $R[T]$ by inverting all the monic polynomials.

\smallskip

\quest
\emph{Let $R$ be a commutative Noetherian ring of dimension $d\geq 2$. Let $I\subset R[T]$ be an ideal 
of height $n$ such that $\mu(I/I^2)=n$. Let $I=(f_1,\cdots,f_n)+I^2$ be given. Assume that  
$IR(T)=(G_1,\cdots,G_n)$ with $G_i-f_i\in I^2 R(T)$. Then,
do there exist $F_1,\cdots,F_n\in I$ such that $I=(F_1,\cdots,F_n)$ with $F_i-f_i\in I^2$?}

\medskip

 \rmk
We first remark on the case $d=n=2$. In this case, the answer to the above question is negative (see \cite[Example 3.15]{brs1}). However, in this case it follows from \cite[Section 7]{d1} that 
there exist $F_1,F_2\in I$ and $\alpha\in SL_{2}(R[T]/I)$ such that $I=(F_1,F_2)$ with
$(\ol{f_1},\ol{f_2})\alpha=(\overline{F_1},\ol{F_2})$, where `bar' denotes images in $I/I^2$.

\smallskip

In view of the above remark, from now on we assume that $3\leq n\leq d\leq 2n-3$
and write this simply as $2n\geq d+3$.
To recast the above question in the language of the Euler class groups let us recall that 
if $R$ is a regular domain (containing an infinite perfect field) of dimension $d\geq 3$, the $n$-th Euler class 
group $E^n(R[T])$ of $R[T]$  has been defined in \cite{drs2} (We remind 
the reader that this definition is not a trivial extension of the definition of the Euler class group
$E^n(R)$ as available in \cite{brs4}).  Let $I\subset R[T]$ be an ideal of height $n$ such that 
$\mu(I/I^2)=n$.
Let $\omega_I:(R[T]/I)^n\sur I/I^2$ be a surjection. One can associate an element $(I,\omega_I)$ in 
$E^n(R[T])$.
The following is the relevant result for us.

\bt\cite[Theorem 3.1]{drs2}
Let $R$ be a regular domain of dimension $d$ containing an infinite perfect field and $n$ be an integer such that 
 $2n\geq d+3$. Let $I\subset R[T]$ be an ideal of height $n$ such that $\mu(I/I^2)=n$.
Let $\omega_I:(R[T]/I)^n\sur I/I^2$ be a surjection. Then $(I,\omega_I)=0$ in $E^n(R[T])$ if and only if there is a surjection $\beta:R[T]^n\sur I$ such that $\beta$ lifts $\omega_I$.
\et

We now consider 
$R(T)$. As the extension 
$R[T]\lra R(T)$ is flat, there is a canonical group homomorphism $\Gamma:E^n(R[T])\lra E^n(R(T))$. The following question is a reformulation of the question posed above.

\medskip

\quest
Is the canonical map $\Gamma:E^n(R[T])\lra E^n(R(T))$ injective?

\smallskip

Here we settle the question in the affirmative when $R$ is a regular domain of dimension $d\geq 3$ which is essentially of finite type over an infinite perfect field $k$ of characteristic $\not =2$. We first prove an easy proposition. For the definitions of $Q'_{2n}(-)$ and $\pi_0(Q'_{2n}(-))$, see Section \ref{homotopy}.

\bp\label{q'}
Let $A$ be any commutative Noetherian ring and $n\geq 2$. Then there is a bijection
$\mu:\pi_0(Q'_{2n}(A))\simeq \pi_0(Q'_{2n}(A[T]))$.
\ep

\proof
For $v\in Q'_{2n}(A)$ we define $\mu([v])=[v]$, sending the orbit $[v]$ to 
$[v]\in \pi_0(Q'_{2n}(A[T]))$. To see that $\mu$ is well defined,  let us first assume that 
$w\in Q'_{2n}(A)$
is homotopic to $v$. Then there is $F(X)\in Q'_{2n}(A[X])$ such that 
$F(0)=v$ and $F(1)=w$. Treat $F(X)$ as an element of $Q'_{2n}(A[T,X])$ and observe
that $[v]=[w]$ in $\pi_0(Q'_{2n}(A[T]))$. 

Now, if $[v]=[w]$ in $\pi_0(Q'_{2n}(A))$, then there
will be a finite sequence of homotopies connecting $v$ and $w$. For each such homotopy we 
do exactly what we did above. Then as elements of $Q'_{2n}(A[T])$, $v$ and $w$ will be connected 
by a finite sequence of homotopies and therefore, $[v]=[w]$ in $\pi_0(Q'_{2n}(A[T]))$.

To prove that $\mu$ is injective, let $v,w\in Q'_{2n}(A)$ be such that 
$[v]=[w]$ in $\pi_0(Q'_{2n}(A[T]))$. Again, it is enough to treat the case when 
$v$ and $w$ are homotopic in $Q'_{2n}(A[T])$ and let us assume so.
Then there is $F(T,X)\in Q'_{2n}(A[T,X])$ such that $F(T,0)=v$ and $F(T,1)=w$.  We then
also have $F(0,0)=v$ and $F(0,1)=w$. Clearly, $G(X)=F(0,X)\in  Q'_{2n}(A[X])$ and
is the desired homotopy. This shows that $\mu$ is injective.

Let $V(T)\in Q'_{2n}(A[T])$. We show that $V(T)$ is homotopic to $V(0)\in Q'_{2n}(A)$.
It is easy to see that $V(TX)\in Q'_{2n}(A[T,X])$. The substitution $X=0$ gives $V(0)$,
whereas the substitution $X=1$ gives $V(T)$. Therefore, $\mu$ is surjective and the proof is complete.
\qed

\smallskip

Let $R$ be a regular domain of dimension $d\geq 3$ which is essentially of finite type over an
infinite perfect  field $k$. Let $n$ be an integer such that  
$2n\geq d+3$.
Recall that in Section \ref{bijection} we  established a bijection 
$\theta_n:E^n(R)\iso \pi_0(Q_{2n}(R))$. If we  assume that $\text{Char}(k)\not = 2$,
then we also have a bijection $\beta_n:\pi_0(Q_{2n}(R))\iso \pi_0(Q'_{2n}(R))$. Therefore,
$\beta_n\theta_n:E^n(R)\iso \pi_0(Q'_{2n}(R))$ is a bijection. From now on we  assume that our rings contain 
$\frac{1}{2}$.

\bt
Let $R$ be a regular domain of dimension $d\geq 3$ which is essentially of finite type over an
infinite perfect  field $k$ of characteristic $\not =2$. Let $n$ be an integer such that  
$2n\geq d+3$. Then we have a bijection $E^n(R[T])\simeq \pi_0(Q'_{2n}(R[T]))$. 
\et

\proof
As $R$ is a regular domain  which is essentially of finite type over an
infinite perfect  field, by \cite[Theorem 3.8]{drs2}  
the group $E^n(R)$ is canonically isomorphic to $E^n(R[T])$. We combine this with
$\beta_n\theta_n:E^n(R)\iso \pi_0(Q'_{2n}(R))$  and use the above proposition.
\qed

\medskip

We now  have the following commutative diagram:
$$
\xymatrix{
         & E^n(R[T]) \ar^{\sim} [r] \ar [d] & \pi_0(Q'_{2n}(R[T])) \ar
       [r]^{\sim} \ar [d] &
        \displaystyle\frac{Q'_{2n}(R[T])}{EO_{2n+1}(R[T])}  \ar [d] &   \\
        & E^n(R(T)) \ar [r]^{\sim} & \pi_0(Q'_{2n}(R(T))) \ar [r]^{\sim} &
        \displaystyle\frac{Q'_{2n}(R(T))}{EO_{2n+1}(R(T))}  & 
        }$$

We can thus transfer the question on the Euler class group to a question on the
map $Q'_{2n}(R[T])/EO_{2n+1}(R[T])\to Q'_{2n}(R(T))/EO_{2n+1}(R(T))$.
 But note that the horizontal bijections are only set-theoretic maps.
However, as we will see below,  for us it will be enough to prove the following statement: 
For $v\in Q'_{2n}(R[T])$,
$$v\stackrel{EO_{2n+1}(R(T))}{\sim} (0,\cdots,0,1)\Rightarrow v\stackrel{EO_{2n+1}(R[T])}
{\sim} (0,\cdots,0,1)$$

We proceed to prove the above assertion in a more general set-up.
 We will need the following factorization lemma for the elementary orthogonal group.

\bl\cite[Lemma 2.3]{st}\label{split} Let $A$ be a commutative Noetherian ring with $2A=A$.
Let  $n\geq 2$ and let $f,g\in A$ be such that 
$fA+gA=A$. If 
$\sigma\in EO_{2n+1}(A_{fg})$,
then there exist $\ga\in  EO_{2n+1}(A_f)$ and $\gb\in EO_{2n+1}(A_g)$ such that
$\sigma=\ga_g\, \gb_f$.
\el

We now prove the following monic inversion principle which is an analogue of \cite[Theorem 1.1]{r2}.
 It will be interesting if one can remove the
regularity assumption from the version we prove.

\bt\label{ravi2}
Let $A$ be a regular domain containing a perfect field $k$ of characteristic $\not = 2$.
Let $n\geq 2$ and $v=(f_1,\cdots,f_n,g_1,\cdots,g_n,h)\in Q'_{2n}(A[T])$. Suppose that
there is a monic polynomial $f\in A[T]$ such that $v\equiv (0,\cdots,0,1)$ mod $EO_{2n+1}(A[T]_f)$.
Then $v\equiv (0,\cdots,0,1)$ mod $EO_{2n+1}(A[T])$.
\et

\proof
Under the assumptions of the theorem, by Theorem \ref{manmi}, $v$ is extended from $A$. More precisely, from Theorem \ref{manmi} we actually have
$$v\stackrel{ O_{2n+1}(A[T],T)}{\sim} v(0),$$
and therefore, by Theorem \ref{stav} it follows that 
$v\stackrel{ EO_{2n+1}(A[T])}{\sim} v(0).$
As $v(0)$ and $v(1)$ are homotopic, applying Theorem \ref{uni1} we have
$$v(0)\stackrel{ EO_{2n+1}(A)}{\sim} v(1), \text{
implying that }
v\stackrel{ EO_{2n+1}(A[T]}{\sim} v(1)\,\,\,\,\,\,\,\,\, (\ast)$$

We first treat the special case when  $f=T$. In this case, 
$v_T\sim (0,\cdots,0,1)$ implies that $v \sim v(1)$ and $v(1) \sim (0,\cdots,0,1)$, and
therefore, we are done:
$$v\stackrel{ EO_{2n+1}(A[T])}{\sim} (0,\cdots,0,1).$$ 

We now show that one can reduce to the special case above. Assume that $f$ is an arbitrary monic polynomial. By ($\ast$) we now have $v=v(1)$ and that
 $v(1) \gs = (0,\cdots,0,1)$ for some $\gs\in EO_{2n+1}(A[T]_f)$. 

Let 
$f^{\ast}=T^{-\text{ deg }f}.f\in A[T^{-1}]$. Then $f^{\ast}(T^{-1}=0)=1$ and note that in $A[T^{-1}]$, we have $f^{\ast}A[T^{-1}]+T^{-1}A[T^{-1}]=A[T^{-1}]$. Further, 
$A[T^{-1},T]_{f^{\ast}}=A[T,T^{-1}]_{f}$.

Applying 
Lemma \ref{split} for the comaximal elements $f^{\ast}$ and $T^{-1}$ in $A[T^{-1}]$,
we have $\gs_T=(\gs_1)_{T^{-1}}(\gs_2)_{f^{\ast}}$, where  
$\gs_1\in EO_{2n+1}(A[T^{-1}]_{f^{\ast}})$ and $\gs_2\in EO_{2n+1}(A[T^{-1},T])$. Now we treat $v(1)$ as an element in 
$Q'_{2n}(A[T^{-1}]_{f^{\ast}})$.  Observe that  the vectors
$v(1)\gs_1$ and $(0,\cdots,0,1)\gs_2^{-1}$ coincide over $A[T^{-1},T]_{f^{\ast}}$. We can
patch them to obtain $w\in Q'_{2n}(A[T^{-1}])$. Note that 
$w_{T^{-1}}=(0,\cdots,0,1)\gs_2^{-1}$ in $A[T,T^{-1}]$. In other words,
$$w_{T^{-1}}\stackrel{ EO_{2n+1}(A[T,T^{-1}])}{\sim} (0,\cdots,0,1)$$
Now apply the special case treated above to conclude that
$w\stackrel{ EO_{2n+1}(A[T^{-1}])}{\sim} (0,\cdots,0,1).$
For this final moment, let us write $Y=T^{-1}$.
As $f^{\ast}(Y=0)=1$ and as $w_{f^{\ast}}=v(1)\gs_1$, specializing at $Y=0$
we obtain 
$$w(0)=v(1)\gs_1(0)\stackrel{EO_{2n+1}(A)}{\sim} (0,\cdots,0,1),$$
showing that $v=v(1)\stackrel{EO_{2n+1}(A)}{\sim} (0,\cdots,0,1)$.
\qed

\bc\label{Mon}
Let $A$ be a regular domain containing a perfect field $k$ of characteristic $\not = 2$.
Let $v\in Q'_{2n}(A[T])$ be such that $v$ is homotopic to $(0,\cdots,0,1)$ in $Q'_{2n}(A(T))$.
Then $v$ is homotopic to $(0,\cdots,0,1)$ in $Q'_{2n}(A[T])$. 
\ec

\proof
From Theorem \ref{uni1} we can conclude that $v\stackrel{EO_{2n+1}(A(T))}{\sim} (0,\cdots,0,1)$.
We can then find a single monic polynomial $f\in A[T]$ such that 
$$v\stackrel{EO_{2n+1}(A[T]_f)}{\sim} (0,\cdots,0,1)$$
Applying the above theorem and Theorem \ref{uni1} again  we are done.
\qed

\medskip

 We are now ready to prove that the monic inversion principle holds for
the Euler class groups.

\bt\label{monicinv}
Let $R$ be a 
regular domain of dimension $d\geq 3$ which is essentially of finite type over an infinite perfect field $k$ of characteristic $\not = 2$. Let $n$ be an integer such that  
$2n\geq d+3$.
Then the natural map $E^n(R[T])\lra E^n(R(T))$ is injective.
\et

\proof
The canonical map $E^n(R[T])\lra E^n(R(T))$ is a group homomorphism. Therefore, 
it is enough to show that if
 $(J,\omega_J)\in E^n(R[T])$ is such that its image $(JR(T),\omega_J\otimes R(T))$ is 
trivial in $E^n(R(T))$, then $(J,\omega_J)=0$ in $E^n(R[T])$. Note that  the horizontal bijections in the following commutative diagram preserve base points.
$$
\xymatrix{
         & E^n(R[T]) \ar^{\sim} [r] \ar [d] & \pi_0(Q'_{2n}(R[T]))  \ar [d]
         &
           &   \\
        & E^n(R(T)) \ar [r]^{\sim} & \pi_0(Q'_{2n}(R(T)))  &
          & 
        }$$
				Therefore we are done by Corollary \ref{Mon}.
\qed

\smallskip

The assertion of Theorem \ref{monicinv}, translated back to  simple ideal-theoretic terms, gives us the main theorem:

\bt\label{monicinv2}
Let $R$ be a regular domain of dimension $d$ which is essentially of finite type over an infinite perfect field $k$ of characteristic $\not = 2$. Let $n$ be a positive integer such that 
$2n\geq d+3$.
Let $I\subset R[T]$ be an ideal of height $n$ such that 
$$I=(f_1,\cdots,f_n)+I^2.$$ 
Assume that  
$IR(T)=(G_1,\cdots,G_n)$ with $G_i-f_i\in I^2R(T)$ for $i=1,\cdots,n$. Then,
there exist $g_1,\cdots,g_n\in I$ such that 
$$I=(g_1,\cdots,g_n), \text{  and } g_i-f_i\in I^2 \text{  for } 
i=1,\cdots,n.$$
\et

\medskip

\noindent
{\bf Acknowledgements.} We sincerely thank Ravi Rao for clarifying numerous queries.  We are deeply indebted to  S. M. Bhatwadekar for his critical reading of an earlier version and for pointing out a mistake in
the proof of Theorem \ref{biject}. The question tackled in this article was proposed to the first
named  author 
by S. M. Bhatwadekar and Raja Sridharan  around the year 2000, as a part of his thesis problem. The first named  author takes this opportunity to thank them once again for their care and encouragement.
The third named author  acknowledges Department of Science and Technology for their INSPIRE research grant.

\end{document}